\documentclass{birkjour}
\usepackage[noadjust]{cite}
\usepackage{xcolor}
\RequirePackage[all]{xy}
\usepackage{amssymb,amsxtra}
\usepackage{enumitem} 
\setlist[enumerate,1]{label=\textup{(\arabic*)}}

\usepackage[lite]{amsrefs}
\usepackage[unicode=true,
colorlinks=true,
linkcolor=blue!75!black,
citecolor=blue!75!black,
urlcolor=purple,
breaklinks=true]{hyperref}

\newtheorem{thm}{Theorem}
\newtheorem{cor}[thm]{Corollary}
\newtheorem{lem}[thm]{Lemma}
\newtheorem{prop}[thm]{Proposition}
\theoremstyle{definition}
\newtheorem{defn}[thm]{Definition}
\theoremstyle{remark}
\newtheorem{rem}[thm]{Remark}
\newtheorem {notat}[thm] {Notation}

\newtheorem*{ex}{Example}

\newcommand{\BibTeX}{B\kern-0.1emi\kern-0.017emb\kern-0.15em\TeX}
\newcommand{\XYpic}{$\mathrm{X\kern-0.3em\raisebox{-0.18em}{Y}}$-$\mathrm{pic}\,$}

\newcommand{\cl}{C \kern -0.1em \ell}  

\mathchardef\mhyphen="2D 
\newcommand*{\UMult}{\mathcal{UM}}
\newcommand*{\Mult}{\mathcal M}

\DeclareMathOperator{\Aut}{Aut}

\DeclareMathOperator{\Iso}{Iso}

\DeclareMathOperator{\PBij}{PBij}
\DeclareMathOperator{\PHomeo}{PHomeo}
\DeclareMathOperator{\PAut}{PAut}
\DeclareMathOperator{\id}{id}
\DeclareMathOperator{\Null}{\mathcal{N}}

\newcommand{\RR}{\mathcal R}
\DeclareMathOperator{\clsp}{\overline{span}}

\newcommand{\ed}{\end{document}}
\begin{document}

%
%
%
%
%
%
%
%
%

\title
 {Banach algebra crossed products by inverse semigroup actions}
 \author[Krzysztof Bardadyn]{Krzysztof Bardadyn }
 \address{%
  Faculty of Mathematics\\ University of Bia\l{}ystok\\
  Cio\l{}kowskiego 1M\\
  15{-}245 Bia\l{}ystok\\
  Poland}
\email{kbardadyn@math.uwb.edu.pl}

\author{Bartosz K. Kwa\'sniewski}
\address{%
  Faculty of Mathematics\\ University of Bia\l{}ystok\\
  Cio\l{}kowskiego 1M\\
  15{-}245 Bia\l{}ystok\\
  Poland}
 \email{b.kwasniewski@uwb.edu.pl}
%
\begin{abstract}
  We give a self-contained and simplified presentation of the theory
  of covariant representations for inverse semigroup actions on Banach
  algebras, which was recently introduced in
  \cite{Bardadyn-Kwasniewski-McKee1} in the twisted case.  The main
  result of this note is a general universal description of the
  associated Banach algebra crossed product, that allows
  disintegration of all representations of the crossed product.  Such
  a disintegration was studied in \cite{Bardadyn-Kwasniewski-McKee1}
  only for actions on spaces.
\end{abstract}
\label{page:firstblob}
\maketitle
\section{Introduction}

Just as the concept of a group was invented to describe  global symmetries of a given object, inverse semigroups  were born  to describe  symmetries that are local in the sense that they `see' only parts of the object \cite{Lawson}.
Inverse semigroup actions on spaces (aka actions by  `partial symmetries') are in a suitable sense equivalent to \'etale groupoids \cite{Paterson, Exel}.
This link is crucial in the study of groupoid algebras as they can be viewed as crossed products for the corresponding inverse semigroup action.
This is quite well explored in the theory of groupoid $C^*$-algebras \cite{Buss-Exel, Sieben97, Paterson, Exel}. In particular, this led, through the notion of a Cartan subalgebra \cite{Renault}, to the proof that every  classifiable simple $C^*$-algebra  has a (twisted) groupoid model \cite{Li}.
Inspired by this,  the authors of \cite{Bardadyn-Kwasniewski-McKee1} recently introduced  groupoid Banach algebras
by first constructing Banach algebra crossed products for inverse semigroup actions and then relating it to various completions of the groupoid convolution algebra. 

In this note we take a closer look at the crossed product itself. For the sake of clarity we ignore twists considered in \cite{Bardadyn-Kwasniewski-McKee1}.
The main result  is a universal description of the Banach algebra crossed product for inverse semigroup actions on approximately unital Banach algebras, which
 gives a bijective correspondence between representations of the crossed product and covariant 
representations of the action (Theorem \ref{thm:universal_description_of_crossed_product}). The hard part here is to find the right notion of a covariant representation and the   proof that every global action `\emph{disintegrates}'  to the representation of the action. In \cite{Bardadyn-Kwasniewski-McKee1} such a  disintegration was proved using a groupoid model for actions on a commutative $C^*$-algebra $C_0(X)$, which are equivalent to actions on a locally compact Hausdorff space $X$.  We give a direct general proof for actions on noncommuatative Banach algebras.

\section{Prologue: group actions and their  crossed products}\label{sect:prologue}
To explain the nuances we start by discussing  group actions, for which the theory is much simpler and  known to experts.
An
\emph{action of discrete group $G$ on a  Banach algebra} $A$  is a group homomorphism $\alpha:G\to \Aut(A)$, where $\Aut(A)$ is a group of isometric automorphisms of $A$.
For technical reasons we will assume that $A$ is \emph{approximately unital}, i.e. $A$ contains a contractive two-sided approximate unit. 
We denote such an action by $\alpha: G\curvearrowright A$.
Then  the \emph{Banach algebra crossed product} is the Banach space
$$
A\rtimes_\alpha G:=\ell^1(G,A)=\{a:G\to A: \|a\|_1:=\sum_{g\in G}\|a(g)\|<\infty \}
$$ 
with the  multiplication given by 
$$
(a*b)(g):=\sum_{h\in G}a(h)\alpha_h(b(h^{-1}g)), \qquad  g\in G, \ a,b\in \ell^1(G,A).
$$
The algebra $A\rtimes_\alpha G$ has the following description using representations. By a \emph{representation of a Banach algebra} in another Banach algebra we mean a contractive algebra homomorphism.
When the target algebra is the Banach algebra $B(E)$ of all bounded operators on a Banach space $E$ we speak of representations on the Banach space $E$.
A \emph{covariant representation of $\alpha$ on a Banach space} $E$ is a pair $(\pi,v)$ 
		where  $\pi:A\to B(E)$ is a  representation, which is (spatially) non-degenerate in the sense that $\overline{\pi(A)E}=E$,   and  $v:G\to \Iso(E)$ is a group homomorphism into the group $\Iso(E)$ of invertible isometries of $E$ such that
\begin{equation}\label{eq:commutation_relation_action}
		v_g\pi(a)=\pi(\alpha_g(a))v_g,\qquad \text{ for all $a\in A$, $g\in G$}.
\end{equation}
The following proposition was proved in \cite[Propositions 2.7, 2.8]{Zeller-Meier} in the  $C^*$-algebraic setting and in \cite[Lemma 2.8]{BK} for $A=C_0(X)$, but the proof works for general approximately unital Banach algebras.
\begin{prop}\label{prop:integration_disintegration}
Let  $\alpha: G\curvearrowright A$ be a group action on a Banach algebra $A$.
 For any covariant representation $(\pi,v)$ of $\alpha$ the formula
 \begin{equation}\label{eq:group_action_representation_integrated}
    \pi\rtimes v(a)=\sum_{g\in G}\pi(a(g))v_g, \qquad a\in A\rtimes_\alpha G, 
    \end{equation}   
where the series is norm convergent,  defines a representation    $\pi\rtimes v$ of  $A\rtimes_\alpha G$,
and  every non-degenerate representation of $A\rtimes_\alpha G$ is of this form.
\end{prop}
The above proposition says that  $A\rtimes_\alpha G$ is a universal Banach algebra for relations given by covariant representations. 
This property determines  $A\rtimes_\alpha G$ uniquely up to an isometric isomorphism. 
To clarify this, let us introduce a notion of a  \emph{covariant representation of $\alpha$  in a Banach algebra $B$}\label{page:covariant_rep}  as  a pair $(\pi,v)$ 
		where  $\pi:A\to B$ is a  representation, which is (algebraically) non-degenerate in the sense that $\overline{\pi(A)B}=B$,   $v:G\to \UMult(B)$ is a group homomorphism into the group of invertible isometries in the multiplier algebra $\Mult(B)$ and
 the commutation relations \eqref{eq:commutation_relation_action} hold. 
		The following can be deduced from Proposition \ref{prop:integration_disintegration} but it also follows from Theorem \ref{thm:universal_description_of_crossed_product}, 
		see  Corollary
		\ref{cor:group_actions}, below.
\begin{thm}
Let  $\alpha: G\curvearrowright A$ be a group action on a Banach algebra $A$. Then the Banach algebra crossed product $A\rtimes_\alpha G$ has the following two properties:
\begin{enumerate}
\item\label{enu:universal_group_of_crossed_product1} there is a covariant representation 
$(\iota, u)$ of $\alpha$ in $A\rtimes_\alpha G$ such that $A\rtimes_\alpha G=\clsp \{\iota(a)u_g: a\in A, g\in G\}$;
\item\label{enu:universal_group_of_crossed_product2} for any covariant representation $(\pi,v)$ of $\alpha$ in a Banach algebra $B$ there is a unique 
representation $\pi\rtimes v:A\rtimes_\alpha G\to B$ such that 
$$
\pi\rtimes v(\iota(a)u_g)=\pi(a)v_g\quad \text{ for all }a\in A, g\in G.
$$
\end{enumerate}
Any Banach algebra satisfying \ref{enu:universal_group_of_crossed_product1} and \ref{enu:universal_group_of_crossed_product2} above is canonically isometrically isomorphic to 
$A\rtimes_\alpha G$.
\end{thm}
	Note that the two notions of non-degeneracy of a representation that we used above differ. Namely, the algebraic non-degeneracy of $\pi:A\to B(E)$ forces $A$ and $\pi$ to be unital, while spatial non-degeneracy is much weaker (e.g. consider the identity embedding of compact operators). 
	We view these non-degeneracy assumptions as \emph{normalization conditions}. The subtleties associated with them become even more troublesome  in the inverse semigroup setting, 
	as there we need to work with ideals, `local multipliers' and their representations.  
In 	particular, we need to decide what a normalization condition should be (and this may depend on the context).  
Symptomatically, already in the $C^*$-algebra setting there are two competing definitions of covariant representations:
one using Hilbert spaces \cite{Sieben97} and one using enveloping $W^*$-algebras \cite{Buss-Exel}. 
Our theory explains how to pass from one representation to the other, and as consequence the resulting algebras are the same (Corollary \ref{cor:Sieben_Buss}).

\section{Covariant representations of inverse semigroup actions}
An \emph{inverse semigroup} is a  semigroup $S$ such that for every $t\in S$ there is a unique `generalized inverse', i.e. $t^*\in S$ such that $t = t t^* t$ and $t^* = t^* t t^*$. 
 This uniqueness is equivalent to that  the idempotents in $S$ commute and so they form a commutative semigroup  $E(S):= \{ e \in S : e^2 = e \}$. A group $G$ is an inverse semigroup with a unique idempotent, which then is necessarily a unit. Idempotents in $S$ induce  a  partial order on $S$ where for any $s,t \in S$ we write $s \le t$ if $s = t s^*s$, which is equivalent 
to  $s = t e$ for some $e \in E(S)$. We refer to \cite{Lawson, Paterson, Exel} for more details. 
\begin{ex}
A \emph{partial bijection} of a set $X$ is a bijection $\theta: U\to \theta(U)$ between two subsets $U,\theta(U)\subseteq X$.
We compose  partial bijections $\theta: U\to\theta(U)$, $\vartheta: V\to \vartheta(V)$ on the largest domain the composition makes sense. Hence   
$$
\theta\circ\vartheta: 
\vartheta^{-1}(U\cap \vartheta(V))\to \theta(U\cap \vartheta(V)).
$$
With this operation  the set $\PBij(X)$  of all partial bijections on $X$ becomes an inverse semigroup,
where   $\theta^*=\theta^{-1}:\theta(U)\to U$ is the inverse to $\theta$, and idempotents  
$E(S)=\{\theta: \theta = id|_U \text{ for some } U\subseteq X  \}$ 
may be identified with subsets of $X$. Moreover,  $\theta\leq\vartheta$ if and only if $\theta$ is the restriction of the map $\vartheta$ to the domain of $\theta$.
This example is generic in the sense that every inverse semigroup $S$ can be identified with an inverse subsemigroup of $\PBij(X)$ for some $X$ (in fact one can take $X=S$ here). 
\end{ex}

By an \emph{ideal} in a Banach algebra we always mean a closed two-sided ideal.
A \emph{partial automorphism} on a Banach algebra $A$ is an isometric algebra isomorphism $\alpha : I \to J$ between two ideals $I,J$ of $A$.  
The set $\PAut(A)$ of partial automorphisms of $A$ together with composition of partial maps
forms an inverse semigroup. 
When $A=C_0(X)$, for a locally compact Hausdorff space $X$, then ideals in $A$ correspond to open sets, and partial isomorphisms
 are given by composition with partial homeomorphisms. In particular, we have
a natural isomorphism of semigroups $\PAut(C_0(X))\cong \PHomeo(X)$ where $\PHomeo(X)$ is the inverse semigroup of partial homeomorphisms.

\begin{defn} 
    An \emph{action of an inverse semigroup $S$ on a Banach algebra $A$} is a  semigroup homomorphism $\alpha: S\to\PAut(A)$  which is `approximately' unital.
		More precisely, $\alpha = \{\alpha_t\}_{t\in S}$ is a family of partial automorphisms $\alpha_t : I_{t^*} \to I_{t}$ of $A$ such that
		\begin{enumerate}[labelindent=40pt,label={(PA\arabic*)},itemindent=1em]
\item\label{enu:inverse_semigroup_action_algebraic1} $		\alpha_s\circ \alpha_t =\alpha_{st}$ (as partial maps) for all $s,t\in S$;
\item\label{enu:inverse_semigroup_action_algebraic2}  each ideal $I_{t}$, $t\in S$, has an approximate unit and the union $\bigcup_{e\in E(S)} I_{e}$ is linearly dense in $A$.  
\end{enumerate}
If $S$ has  zero, then we will also assume that $I_0=\{0\}$.
We will denote such inverse semigroup actions by writing 	$\alpha: S\curvearrowright A$.
\end{defn}
\begin{rem}
If $S$ is unital, then in the presence of \ref{enu:inverse_semigroup_action_algebraic1}  the second part of \ref{enu:inverse_semigroup_action_algebraic2}  is equivalent to assuming that  $\alpha$ is unital, i.e. $\alpha_1 = \id_A$.
\end{rem}
\begin{rem}\label{rem:twisted_relations} When $A$ is a $C^*$-algebra the above definition appeared first in the work of Sieben~\cite[Definition 3.1]{Sieben97} (every ideal in a $C^*$-algebra automatically has an approximate unit). For twisted actions on Banach algebras it was introduced in \cite{Bardadyn-Kwasniewski-McKee1}.
Preservation of the semigroup law required in \ref{enu:inverse_semigroup_action_algebraic1} has a number of consequences.
In particular,  for any action $\alpha: S\curvearrowright A$,  for all $s,t\in S$ and $e\in E(S)$ we have
\[
  \alpha_s (I_{s^*}\cap I_t) = I_{st},\quad   I_t = I_{tt^*},  \quad \alpha_e = \text{id}|_{I_e}, \quad \alpha_{t^*}=\alpha^{-1}_t,
\] 
and $I_{s}\subseteq I_{t}$ whenever $s\leq t$. 
We will use these relations without warning.
\end{rem}
\begin{ex}
Inverse semigroup actions $\alpha$ on  $A=C_0(X)$, where $X$ is a locally compact Hausdorff space, correspond to
inverse semigroup actions $\theta$ on $X$ as defined in \cite{Exel}, where  $\alpha_{t}:C_0(X_{t^*})\to C_0(X_{t})$ is given by $\alpha_t(a)=a\circ \theta_{t^*}$ for $a\in C_0(X_{t^*})$, and 
$\theta_{t*}:X_{t}\to X_{t^*}$ is partial homeomorphism of $X$.
\end{ex}
\begin{ex}
If $S=G$ is a group and $A$ is a Banach algebra, then inverse semigroup actions of $S$ on $A$ coincide with group actions of $G$ on $A$.
\end{ex}
The following definition is a straightforward generalization of the corresponding $C^*$-algebraic notion ~\cite[Definition 3.4]{Sieben97}, as well as 
the Banach algebra definition for group actions we used in Section \ref{sect:prologue}. 
For any Banach algebra $B$ we denote by $B_1:=\{b\in B: \|b\|\leq 1\}$ the closed unit ball.
\begin{defn}\label{defn:covariant_representations_on_Banach_spaces}
    Let $\alpha: S\curvearrowright A$ be an inverse semigroup action. A \emph{covariant representation of $\alpha$ on a Banach space $E$} is a pair $(\pi, v)$ where $\pi:A\to B(E)$ is a representation and $v: S\to B(E)_1$ is a semigroup homomorphism into contractive operators such that
    \begin{enumerate}[labelindent=40pt,label={(SCR\arabic*)},itemindent=3em]
        \item\label{enu:covariant_representations_on_Banach_spaces1} $v_t\pi(a)=\pi(\alpha_t(a)))v_{t}$ for all $a\in I_{t^*}$,  $t\in S$;
        \item\label{enu:covariant_representations_on_Banach_spaces2} $v_tE=\overline{\pi(I_{t})E}$ for $t\in S$.
    \end{enumerate}
\end{defn}
\begin{rem}\label{rem:covariant_representations_on_Banach_spaces}
 Since $I_{t}=I_{tt^*}$ for all $t\in S$ and all idempotents in $S$ are of the form $tt^*$,   we see that   \ref{enu:covariant_representations_on_Banach_spaces2} is equivalent to assuming that 
the range of the contractive projection $v_e$ corresponding to an idempotent $e\in E(S)$ is  $\overline{\pi(I_{e})E}$. 
This condition together with \ref{enu:covariant_representations_on_Banach_spaces1} and the fact that $\alpha_e=\id_{I_e}$, for  all $e\in E(S)$, imply that  
$
\pi(a)=v_{e}\pi(a)=\pi(a)v_{e}
 \text{ for  all $a\in I_{e}$, $e\in E(S)$.}
$
\end{rem}
The operators $\{v_t\}_{t\in S}$ in the above definition are partial isometries in the sense of Mbekhta \cite{Mbekhta}, as they are contractive and 
$$
v_t v_{t^*} v_{t}=v_{t}, \qquad  v_{t^*} v_{t} v_{t^*}=v_{t^*}, \qquad t\in S. 
$$
In particular, for representations on Hilbert space we necessarily have $v_{t^*}=v_t^*$.
For actions on $C^*$-algebras or more generally when ideals that appear in the action have hermitian approximate units,  operators $\{v_t\}_{t\in S}$ are necessarily Moore-Penrose partial isometries,
cf.  \cite[Corollary 3.20]{Bardadyn-Kwasniewski-McKee1}. 

Construction of the $C^*$-crossed product  for an inverse semigroup action  \cite{Sieben97} can be readily adopted to produce a  Banach algebra which is universal in the sense that every 
covariant representation $(\pi,v)$  of the action integrates  to a representation $\pi\rtimes v$ of the algebra. 
But when it comes to disintegration the situation becomes much more subtle.
In general, we will use the \emph{double dual} $B''$ of a Banach algebra $B$. 
The reason is that `local multipliers' live in $B''$ rather than in $B$ or even in $\Mult(B)$. Also we will need some sort of weak topology and  $B''$   has a weak$^*$ topology induced by $B'$, which we will refer to as $B'$-topology. 
Recall \cite{Dales} that $B''$ is again naturally a Banach algebra with either of the Arens products, in which $B$ sits as a $B'$-weakly dense Banach algebra.
Namely, for $a,b \in B''$, the \emph{first Arens product} is defined as
    $$
    a \square b=B^\prime\mhyphen\lim_{\alpha}B^\prime\mhyphen\lim_{\beta}a_\alpha b_\beta
    $$
    where $\{a_\alpha\}_{\alpha}$, $\{b_\beta\}_{\beta}$ are nets $B^\prime$-convergent to $a$, $b$ respectively. Similarly, the \emph{second Arens product} is given by
    $$
    a\cdot b=B^\prime\mhyphen\lim_{\beta}B^\prime\mhyphen\lim_{\alpha}a_\alpha b_\beta.
    $$
    The Banach space $B^{\prime\prime}$ with any of the above products is a Banach algebra. 
			In general these products are different. If the two products coincide, then  $B$ is called \emph{Arens regular}, 
see \cite{Dales} for more details. 
For the sake of clarity, in this article \emph{we will only use the second Arens product}.
Our decision is arbitrary, and in any case we will be primarily interested in products $a v\in B$ where $a\in B$ and $v\in B''$, 
in which case the first and second Arens  products always agree. 
Also recall that every $C^*$-algebra is Arens regular. 

For any representation $\pi : A \to B$ the double adjoint $\pi'' : A'' \to B''$ is a  $B'$-weakly continuous representation that extends $\pi$, and $\pi''$ is isometric whenever $\pi$ is.
In particular, for every Banach subalgebra $A\subseteq B$ we may identify $A''$ with a Banach subalgebra of $B''$.
Also any   approximate unit  $\{\mu_{i}\}$ in $A$ converges $A'$-weakly  to a left identity $1_A$ in $A''$ (and a right identity for the first Arens product in $A''$) which is a unit in the multiplier algebra $\Mult(A):=\{ b \in A'' : \text{$ba, ab \in A$ for all $a\in A$} \} \subseteq A''$, see \cite[Proposition 2.9.16 and Theorem 2.9.49]{Dales}.
In particular, for each ideal $I$ in $A$, which has an approximate unit, we will adopt the identifications $I\subseteq \Mult(I)\subseteq I''\subseteq A''$.

\begin{defn}\label{defn:covariant_representation_in_Banach_algebras}
    Let $\alpha: S\curvearrowright A$ be an inverse semigroup action. A \emph{covariant representation of $\alpha$ in a Banach algebra $B$} is a pair $(\pi,v)$ where $\pi: A\to B$ is a representation and $v: S\to (B^{\prime\prime})_1$  is a map satisfying
    \begin{enumerate}[labelindent=40pt,label={(CR\arabic*)},itemindent=2em]
        \item\label{item:covariant_representation1} $v_t\pi(a)= \pi(\alpha_t(a))v_t\in B$ for all $a\in I_{t^*}$;
        \item\label{item:covariant_representation2} $\pi(a)v_sv_t=\pi(a)v_{st}$ for all $a\in I_{st}$, $s$, $t\in S$;
        \item\label{item:covariant_representation3} $\pi(a)v_e=\pi(a)$ for all $a\in I_e$, $e\in E(S)$.
    \end{enumerate}
  Such a covariant representation $(\pi, v)$ is non-degenerate, if $\overline{\pi(A)B}=B$. We call $B(\pi,v):=\clsp\{ \pi(a)v_t: a\in I_{t}, t\in S\}$ the \emph{range of $(\pi,v)$}.
\end{defn}
\begin{rem}\label{rem:general_covariant_rep}
Condition~\ref{item:covariant_representation1} is equivalent to $\pi(a)v_t=v_t\pi(\alpha_{t^*}(a))\in B$ for all $a\in I_{t}$, $t\in S$. 
Note that this condition means two things. Firstly, this is a commutation relation. Secondly it requires that the considered product is in $B$, rather than in $B''$. In particular it implies that  $B(\pi,v)\subseteq B$.
Conditions~\ref{item:covariant_representation1} and \ref{item:covariant_representation3} imply that $v_{e}\pi(a)=\pi(a)$ for all $a\in I_{e}$, $e\in E(S)$ (note that $\alpha_e=\id_{I_e}$). Since $I_t=I_{tt^*}$, \ref{item:covariant_representation2} and \ref{item:covariant_representation3} imply that $\pi(a)v_t^*v_t=\pi(a)v_{t^*t}=\pi(a)$ for all $a\in I_t$. Employing also \ref{item:covariant_representation1} we get 
\begin{equation}\label{eq:CovarianceConditionAlternate}
    v_t\pi(a)v_{t^*} = \pi \big( \alpha_t(a) \big) \qquad \text{for all $a\in I_{t^*}$, $t\in S$.} 
\end{equation}
\end{rem}

\begin{lem}\label{lem:range_of_covariant_rep}
Let $(\pi,v)$ be a covariant representation of $\alpha$ in a Banach algebra $B$.
The range $B(\pi,v)$ is a Banach subalgebra of $B$. 
The spaces  $A_t = \{ \pi(a_t)v_t: a_t\in I_{t} \}$, $t\in S$, form a grading of $B(\pi,v)$ over the inverse semigroup $S$ in the sense that 
\[
    B(\pi,v) = \clsp \{ A_t : t\in S \}, \qquad A_s A_t \subseteq A_{st}, \qquad \text{$s\leq t$ implies $A_s\subseteq A_t$} ,
\]
for all $s,t \in S$. 
In fact, for all $a_t\in I_{t}, a_s\in I_{s}$, $s,t\in S$, we have
\begin{enumerate}
    \item\label{enu:range_of_covariant_rep1} $\pi(a_s)v_s \cdot \pi(a_t)v_t = \pi \big( \alpha_s(\alpha_{s^*}(a_s)a_t)\big)v_{st}$ and $\alpha_s \big( \alpha_{s^*} (a_s) a_t \big) \in I_{st}$;
    \item\label{enu:range_of_covariant_rep2} $s\leq t$ implies $I_{s}\subseteq I_{t}$ and $\pi(a_s)v_s = \pi ( a_s ) v_t$.
\end{enumerate}
\end{lem}
\begin{proof} 
\ref{enu:range_of_covariant_rep1}. For any 
$a_t\in I_{t}, a_s\in I_{s}$, $s,t\in S$, we have $\alpha_{s^*}(a_s) a_t\in I_{s^*} I_{t}= I_{s^*}\cap I_{t}$ and so
$\alpha_{s}(\alpha_{s^{*}}(a_s)a_t)\in I_{st}$ because $\alpha_s(I_{s^*}\cap I_{t})=I_{st}$, see Remark \ref{rem:twisted_relations}.
Using \ref{item:covariant_representation1} we get
\begin{align*}
    \pi(a_s)v_s \cdot \pi(a_t)v_t &= v_s \pi \big( \alpha_{s^*}(a_s) \big) \cdot \pi(a_t) v_t = v_s \pi \big( \alpha_{s^*}(a_s) a_t \big) v_t 
		\\
		&= \pi \Big( \alpha_s \big( \alpha_{s}^{*}(a_s) a_t \big) \Big) v_{st}. 
\end{align*}

\ref{enu:range_of_covariant_rep2}. If $s\leq t$, then $I_{s}\subseteq I_{t}$ by Remark~\ref{rem:twisted_relations}. 
By \ref{item:covariant_representation3} and \ref{item:covariant_representation2}, for  $a \in I_{s}$ we get  
$\pi(a)v_t=\pi(a)v_{ss^*} v_t= \pi(a)v_{ss^*t} = \pi(a)v_{s}.
$ 
\end{proof}

Notice that the elements $v_t$ of Definition~\ref{defn:covariant_representation_in_Banach_algebras} are not required to be partial isometries, 
and in particular in general they do not form a semigroup. However, this can always be arranged by `normalizing'. 
Depending on the structure of the target algebra $B$, this normalization can be done in different ways. In general, knowing nothing about $B$,
 we can always use the $B'$-topology of $B''$ as follows.
\begin{defn}
  Let $\alpha: S\curvearrowright A$ be an inverse semigroup action.  A covariant representation $(\pi,v)$ of $\alpha$ in a Banach algebra $B$ is 
	\emph{normalized}, if 
\begin{equation}\label{eq:covariant_representation4}
 v_{t}= B'\mhyphen\lim_{i} (\pi(\mu_{i}^t)v_{t}) \text{  for an an approximate unit   $\{\mu^t_i\}_{i}$ in $I_t$, $t\in S$.}
\end{equation}
	\end{defn}
	\begin{rem}\label{rem_normalized_rep}
Since the elements 	$\pi(\mu_{i}^t)v_{t}$ appearing in \eqref{eq:covariant_representation4} are in the unit ball of $B(\pi,v)\subseteq B(\pi,v)''\subseteq B''$, 
Banach-Alaoglu Theorem implies that if the limit in \eqref{eq:covariant_representation4} exists it sits in $B(\pi,v)''$. In particular, \eqref{eq:covariant_representation4}
is equivalent to 
$$
v_{t}= B(\pi,v)'\mhyphen\lim_{i} (\pi(\mu_{i}^t)v_{t})\, \text{  where $\{\mu^t_i\}_{i}$ is an approximate unit in $I_t$,  $t\in S$.}
$$
	\end{rem}
	\begin{notat}
	We  denote by $1_t$ the unit in $\Mult(I_t)\subseteq I_t''\subseteq A''$. Note that $1_{t}=1_{tt^*}$ and $1_{t^*}=1_{t^*t}$ for all $t\in S$, because $I_{t}=I_{tt^*}$
	and    $I_{t^*}=I_{t^*t}$.
\end{notat}	
\begin{prop}\label{prop:normalized_twisted_rep}
A pair $(\pi,v)$ is a $B'$-normalized covariant representation of $\alpha$ if and only if $\pi : A \to B$ is a representation and $v : S\to (B'')_{1}$ is a semigroup homomorphism   such that $\pi(a) v_t \in B$, for $a\in I_t$, $t\in S$, and
\begin{enumerate}
    \item\label{item:normalized_covariant_representation1} $v_t \pi(a) v_{t^*} = \pi(\alpha_{t}(a))$ for all $a\in I_{t^*}$, $t\in S$;
    \item\label{item:normalized_covariant_representation3} $v_{e} = \pi''(1_{e})$ for all $e\in E(S)$.
\end{enumerate}
Moreover, for any covariant representation $(\pi,v)$ of $\alpha$ in $B$ putting
$$
\tilde{v}_{t} = B'\mhyphen\lim_{i} (\pi(\mu_{i}^t)v_{t}) = v_t \pi''(1_{t^*}),\qquad  t\in S,
$$ 
the pair $(\pi,\tilde{v})$ is  a unique $B'$-normalized covariant representation  such that $\pi(a)v_t=\pi(a)\tilde{v}_{t}$ for all $a\in I_{t}$ and $t\in S$.
In particular,
$
B(\pi,v)=B(\pi,\tilde{v}).
$ 
\end{prop}
\begin{proof} 
Assume $(\pi,v)$ is a pair with the properties described in the assertion. 
For every $t\in S$ and $a\in I_{t^*}$ we have $\pi(a) \pi''(1_{t^*})=\pi(a)$, hence  using also \ref{item:normalized_covariant_representation3} and \ref{item:normalized_covariant_representation1} we have 
\[
    v_t \pi(a) = v_t \pi(a) \pi''(1_{t^*}) = v_t \pi(a) v_{t}^*v_t = \pi \big( \alpha_t(a) \big) v_t  ,
\] 
which is \ref{item:covariant_representation1}. We assume that $v$ is a semigroup homomorphism, which is clearly stronger than \ref{item:covariant_representation2}.
Using the semigroup law and \ref{item:normalized_covariant_representation3} we get 
$$
v_t =v_{tt^*t}= v_t v_{t^*}v_t = \pi''(1_{t}) v_t,
$$ 
which gives \eqref{eq:covariant_representation4}.  Applying  the semigroup law and \ref{item:normalized_covariant_representation3} to an idempotent $e\in E(S)$  we get 
$
v_{e}=v_{ee}=v_{e} v_e=  \pi''(1_e),
$
which implies \ref{item:covariant_representation3}.
Hence $(\pi,v)$ is a $B'$-normalized covariant representation of $\alpha$.

Let now $(\pi,v)$ be any $B'$-normalized covariant representation. 
For each $t\in S$ choose an approximate unit $\{\mu^t_i\}_{i}$ in $I_t = I_{tt^*}$. 
For $a \in I_{t^*}$ we get 
\[
    v_t \pi(a)  v_{t^*} = \pi \big( \alpha_t(a) \big) v_{t}v_{t^*}=
         \pi ( \alpha_t(a) ) v_{tt^*}= \pi( \alpha_t(a) ),
\]
which is \ref{item:normalized_covariant_representation1}. 
 If $e\in E(S)$, then by normalization \eqref{eq:covariant_representation4} and  \ref{item:covariant_representation3} we get
	 $$
	v_e=B'\mhyphen\lim_{i} (\pi(\mu^e_i) v_e)=B'\mhyphen\lim_{i} \pi(\mu^e_i)=\pi''(1_{e}),
	$$
	which is \ref{item:normalized_covariant_representation3}. 
Now take any $s,t\in S$. Recall that $\alpha_s(I_{s^*}\cap I_t)=I_{st}$ and therefore $\{\alpha_s(\alpha_{s^*}(\mu^{s}_j)\mu^{t}_i)\}_{i,j}$ is 
and approximate unit in $I_{st}$. Using this, normalization \eqref{eq:covariant_representation4} and Lemma \ref{lem:range_of_covariant_rep}\ref{enu:range_of_covariant_rep1} we get
\[
\begin{split}
    v_{s} v_{t} &=  v_{s} B'\mhyphen\lim_{i}[ \pi (\mu^{t}_i) v_{t}]  = B'\mhyphen\lim_{i} [v_{s} \pi (\mu^{t}_i) v_{t}] 
        \\
				&= B'\mhyphen\lim_{i} B'\mhyphen\lim_{j} [\pi (\mu^{s}_j) v_{s} \pi (\mu^{t}_i) v_{t}] 
				\\
				&=B'\mhyphen\lim_{i} B'\mhyphen\lim_{j} [\pi \big( \alpha_s(\alpha_{s^*}(\mu^{s}_j)\mu^{t}_i)\big) v_{st}] =v_{st}.
\end{split}
\]
Hence $(\pi,v)$ satisfies all properties in the assertion.

Now let $(\pi,v)$ be any covariant representation of $(\alpha,u)$ in $B$. 
Using \ref{item:covariant_representation1}   we get
\[
\begin{split}
    B' \mhyphen\lim_{i} \big( \pi(\mu^t_i) v_t \big)  &= B'\mhyphen\lim_{i} \big( v_t \pi(\alpha_{t^*}(\mu^t_i)) \big) 
		= v_t B'\mhyphen\lim_{i} \pi \big( \alpha_{t^*}(\mu^t_i) \big) 
		\\
		&=v_{t}\pi''(1_{t^*}),
\end{split}
\]
because $\{\alpha_{t^*}(\mu^t_i)\}_{i}$ is an approximate unit in $I_{t^*}$ which is $B'$-convergent to $1_{t^*}$. 
Thus we that the limit $\tilde{v}_{t} := B'\mhyphen\lim_{i} ( \pi ( \mu_{i}^t ) v_{t} ) = v_t \pi''( 1_{t^*} )$ exists and does not depend on the choice of an approximate unit $\{\mu^t_i\}_{i}$ in $I_t$, $t\in S$. 
Clearly if $(\pi,\tilde{v})$ is a $B'$-normalized covariant representation satisfying $\pi(a) v_t = \pi(a) \tilde{v}_{t}$ for all $a\in I_{t}$ and $t\in S$, then each $\tilde{v}_{t}$ must be given by such a limit.
For $a\in I_{t^*}$ we get  
\[
\begin{split}
    \tilde{v}_t \pi(a) &= v_t \pi''(1_{t^*}) \pi(a)  = v_t\pi(a) = \pi \big( \alpha_{t}(a) \big) v_t
		\\
		&= \lim_{i} \pi \big( \alpha_{t}(a) \big) \pi(\mu_i^{t}) v_{t} = \pi \big( \alpha_{t}(a) \big) \tilde{v}_t.
\end{split}
\]
This shows that $(\pi,\tilde{v})$ satisfies \ref{item:covariant_representation1} and  $\pi(a)v_t=\pi(a)\tilde{v}_{t}$ for all $a\in I_{t}$ and $t\in S$.
Using this we get $\pi(a)\tilde{v}_e =\pi(a)v_e = \pi(a)$ for all $a\in I_{e}$, $e\in E(S)$, so \ref{item:covariant_representation3} holds. 
Finally, for all $s,t\in S$ and $a\in I_{st}$ the calculation 
\[
    \pi(a) \tilde{v}_s \tilde{v}_t = v_s \pi \big( \alpha_{s^*}(a) \big) \tilde{v}_t = \pi(a)v_s v_t = \pi( a) v_{st} = \pi ( a ) \tilde{v}_{st} 
\]
gives \ref{item:covariant_representation2}. As $\tilde{v}$ satisfies  \eqref{eq:covariant_representation4} by construction,
conclude  that $(\pi,\tilde{v})$ is the desired $B'$-normalized covariant representation of $\alpha$.
\end{proof}\begin{rem}\label{rem:normalization_asymetry}
The above proposition implies that for any covariant representation $(\pi,v)$, the normalization  condition \eqref{eq:covariant_representation4}
	is equivalent to assuming that $v_t=v_t\pi''(1_{t^*})$, for all $t\in S$, and it implies that $v_t=\pi''(1_t)v_t$, for all $t\in S$. When $B$ is Arens regular, 
	the latter implication can be reversed. In general it is not clear. This asymmetry results from our arbitrary choice to use the second  Arens'  product in $B''$.
\end{rem}

\begin{cor}\label{cor:range_of_covariant_rep}
Retain the assumptions and notation from Lemma \ref{lem:range_of_covariant_rep}.
Assume in addition that $A$ and $B$ are $C^*$-algebras.
Then 
$$(\pi(a_t) v_t)^* = \pi \big( \alpha_{t^*}(a_{t}^*) \big)v_{t^*} \qquad \text{ for all $a_t\in I_{t}$, $t\in S$,}
$$ 
so $A_t^* = A_{t^*}$, and $\{A_t\}_{t\in S}$ is a saturated grading of $B(\pi,v)$ over  $S$ 
in the sense of \cite[Definition 6.15]{Kwa-Meyer0}, see also \cite[Definition 7.1]{Exel}.
\end{cor}
\begin{proof}

If $B$ is a $C^*$-algebra then $B''$ is again a $C^*$-algebra -- the enveloping $W^*$-algebra of $B$.
Contractive homomorphism $\pi:A\to B$ between $C^*$-algebras is necessarily $*$-preserving.
By passing to the normalized covariant representation $(\pi,\tilde{v})$ described in Proposition~\ref{prop:normalized_twisted_rep} we see that for every $t\in S$, the operator  $\tilde{v}_{t}:=\pi''(1_{t})v_t = v_t\pi''(1_{t^*})$ is a partial isometry whose adjoint in $B''$ is $\tilde{v}_{t^*}=\pi''(1_{t^*})v_{t^*} = v_{t^*}\pi''(1_{t})$. Thus for $a_t\in I_{t}$, $t\in S$, we get
$$
 \big( \pi(a_t) v_t \big)^*=\big( \pi(a_t) \tilde{v_t} \big)^* = \tilde{v}_{t^*} \pi(a_t^*) =
		 \pi(\alpha_{t^*}(a_t^*))\tilde{v}_{t^*}=\pi(\alpha_{t^*}(a_t^*))v_{t^*}.  \qedhere
$$
\end{proof}

\begin{rem}\label{rem:unital_actions}
If  each $I_t$, $t\in S$, is unital, so that $1_t\in A$, $t\in S$, then for any covariant representation $(\pi,v)$  of $\alpha$ in $B$ its $B'$-normalization  $\tilde{v}_t = \pi(1_t)v_t$, $t\in S$, takes values in $B$ by \ref{item:covariant_representation1}. Hence, in this case without changing the range of covariant representations, one can avoid  talking about $B''$ altogether.
In this case  Proposition~\ref{prop:normalized_twisted_rep} 
could be formulated without the use of the bidual algebra $B''$ and the extended representation $\pi''$. 
Namely, when each $I_t$, $t\in S$, is unital,  a normalized covariant representation of $\alpha$ 
is a  pair $(\pi,v)$  where  $\pi : A \to B$ is a representation and $v : S\to B_{1}$ is a semigroup homomorphism   such that 
\begin{enumerate}
    \item\label{item:unital_normalized_covariant_representation1} $v_t \pi(a) v_{t^*} = \pi(\alpha_{t}(a))$ for all $a\in I_{t^*}$, $t\in S$;
    \item\label{item:unital_normalized_covariant_representation3} $v_{e} = \pi(1_{e})$ for all $e\in E(S)$.
\end{enumerate}
\end{rem}
When the inverse  semigroup action $\alpha: S\curvearrowright A$ is arbitrary we    may restrict to maps taking values in $B$ (rather then in $B''$) whenever $B$ has a  predual Banach space $B_*$.  Then we have the canonical embedding $B_*\subseteq  B_*''=B'$ (inclusion is proper if $B$ is not reflexive), and so 
the $B_*$-topology on $B''$ is weaker than $B'$-topology. 
However these topologies coincide on bounded subsets of $B$:

\begin{lem}
A bounded net of elements in $B$ is $B_*$-convergent if and only if it is $B'$-convergent, and then 
the limit is in $B$ rather in $B''$.
\end{lem}
\begin{proof} Let $\{b_i\}_i$ be a bounded net of elements in $B$. Since $B_*\subseteq B'$, if $\{b_i\}_i$  is $B'$-convergent to $b\in B''$, then $\{b_i\}_i$ is  also $B_{*}$-convergent to $b$,
and then Banach-Alaoglu Theorem  implies that $b\in B$ (as the net $\{b_i\}_i$ is contained in a closed ball in $B$). 
Conversely, assume that $\{b_i\}_i$ is $B_{*}$-convergent to $b$ (necessarily in $B$). 
As the net $\{b_i\}_i$ is contained in a closed ball in $B''$, Banach-Alaoglu Theorem  implies that every subnet of $\{b_i\}_i$ contains 
a  $B'$-convergent subnet, which has to converge to $b$ (as $B_*$-topology is weaker than $B'$-topology).
This implies that $\{b_i\}_i$ is $B'$-convergent to $b$. 
\end{proof}
\begin{cor}\label{cor:predual_Banach_space}
If a Banach algebra $B$ has a Banach space predual $B_*$, then a  covariant representation $(\pi,v)$ of $\alpha$ in $B$ is normalized 
if and only if
\begin{equation}\label{eq:covariant_representation5} 
v_{t}= B_*\mhyphen\lim_{i} (\pi(\mu_{i}^t)v_{t}), \text{ for an approximate unit $\{\mu^t_i\}_{i}$   in $I_t$,  $t\in S$,}
\end{equation}
 and then we necessarily have  $\{v_{t}\}_{t\in S}\subseteq B$ (there is no need to consider $B''$).
\end{cor}
Description of normalized covariant representations simplifies even further when  $B$ is a dual Banach algebra. 
Recall that a \emph{dual Banach algebra} is a pair $(B, B_*)$ where $B$ is a Banach algebra,  $B_*$ is a predual Banach space of $B$ and  multiplication in $B$ is $B_*$-weak separately continuous, see \cite{Runde, Daws}.
 Examples of dual Banach algebras include all $W^*$-algebras  with their unique preduals, and all algebras $B(E)$ of all bounded operators on a reflexive Banach space $E$, equipped with the canonical predual Banach space $E' \widehat{\otimes} E$. We use the following result, which is a consequence of \cite[Theorem 5.6]{Ilie_Stokke}.
\begin{prop}[Ilie-Stokke]\label{prop:Ilie-Stokke}
Let $A$ be a Banach algebra with an approximate unit.  Every representation $\pi : A \to B$ into a dual Banach algebra $(B, B_*)$ has a unique extension to a representation $\overline{\pi} : \Mult(A) \to B$ which is strictly-$B_*$-continuous; it is given by  $\overline{\pi}(m) := B_*\mhyphen\lim_{i} \pi(m \mu_{i})$ for any $m$ in $\Mult(A)$ and any approximate unit $\{\mu_{i}\}_i$ in $A$. 
\end{prop}
 
\begin{cor}\label{cor:B*_normalized_twisted_rep} 
 Let $(B,B_*)$ be a dual Banach algebra. 
A pair $(\pi,v)$ is a normalized covariant representation  of an action $\alpha: S\curvearrowright A$ in $B$ if and only if $\pi : A \to B$ is a representation and $v : S\to (B)_{1}$ is a semigroup homomorphism   such that 
\begin{enumerate}
    \item\label{item:B*normalized_covariant_representation1} $v_t \pi(a) v_{t^*} = \pi(\alpha_{t}(a))$ for all $a\in I_{t^*}$, $t\in S$;
    \item\label{item:B*normalized_covariant_representation3} $v_{e} = \overline{\pi}_e(1_{e})$ for all $e\in E(S)$  where $\overline{\pi}_e : \Mult(I_e) \to B$ is the unique strictly-$B_*$-continuous extension of $\pi|_{I_e}$, given by Proposition \ref{prop:Ilie-Stokke}.
\end{enumerate}
Moreover, for any covariant representation $(\pi,v)$ of $\alpha$ in $B$ putting
$$
\tilde{v}_{t} =B_*\mhyphen\lim_{i} \pi(\mu_{i}^t)v_{t} =  \overline{\pi}_{tt^*}(1_{t^*t}) v_t =v_t\overline{\pi}_{t^*t}(1_{t^*t}) ,\qquad  t\in S,
$$ 
the pair $(\pi,\tilde{v})$ is  a unique $B_*$-normalized covariant representation  such that $\pi(a)v_t=\pi(a)\tilde{v}_{t}$ for all $a\in I_{t}$ and $t\in S$.
In particular,
$
B(\pi,v)=B(\pi,\tilde{v}).
$ 
\end{cor}
\begin{proof} Translate Proposition~\ref{prop:normalized_twisted_rep} using Corollary \ref{cor:predual_Banach_space} and Proposition \ref{prop:Ilie-Stokke}.
\end{proof}
Now coming back to covariant representations on a Banach space $E$, as described in Definition \ref{defn:covariant_representations_on_Banach_spaces},  we notice that they can be viewed as covariant representations in  the algebra $B(E)$ which are `normalized' using strong topology.
But if $E$ is reflexive, then they in fact coincide with ($B_*$-)normalized covariant representation in $B(E)$ where $B_*$ is the canonical predual  $B(E)$.

\begin{lem}\label{lem:covariant_reps_implies_MP^partial_isos}
Let $(\pi,v)$ be a covariant representation of  $\alpha: S\curvearrowright A$ on a Banach space $E$. 
Then $(\pi,v)$ is a covariant representation in the Banach algebra $B(E)$ such that each $v_e\in B(E)$, $e\in E(S)$, is a  strong limit of $\{\pi(\mu_i^e)\}_{i}$, where  $\{\mu_{i}^e\}_{i}$ is a contractive approximate unit in $I_{e}$.
\end{lem}
\begin{proof} It is clear that $(\pi,v)$ is a covariant representation in the Banach algebra $B(E)$, see Remark \ref{rem:covariant_representations_on_Banach_spaces}. 
Note that $v_e\in B(E)$, for $e\in E(S)$, is a projection onto $\overline{\pi(I_e)E}$ which commutes with elements of $\pi(I_e)$.  This implies that $v_e$ is a  strong limit of $\{\pi(\mu_i^e)\}_{i}$.
\end{proof}
Recall \cite{Ryan} that if $E$ is a reflexive Banach space, then the projective tensor product Banach space $E' \mathbin{\widehat{\otimes}} E$ is naturally a predual of $B(E)$ making it a dual Banach algebra. Here we identify a simple tensor $\xi\otimes f\in E\widehat{\otimes}E^\prime$ with the functional $\widehat{\xi\otimes f}\in B(E)^\prime$ given by $\widehat{\xi\otimes f}(T):=f(T\xi)$, 
$T\in B(E)$.

\begin{prop}\label{prop:reflexive_covariant_representations}
Assume that $E$ is a reflexive Banach space.
Then normalized covariant representations  of $\alpha$  in the Banach algebra $B(E)$ coincide with  covariant representations of $\alpha$  on the Banach space $E$.
\end{prop}
\begin{proof} We treat  $B(E)$ as a dual Banach algebra with  $B_* := E' \mathbin{\widehat{\otimes}} E$.
Assume that $(\pi,v)$ is a normalized covariant representation in $B(E)$, and so $\{v_t\}_{t\in S}\subseteq B(E)$ are partial isometries with the associated generalized inverses $\{v_t^*\}_{t\in S}\subseteq B(E)$, satisfying the relations described in Corollary~\ref{cor:B*_normalized_twisted_rep}. Since it is a covariant representation it satisfies  \ref{enu:covariant_representations_on_Banach_spaces1} which in this case is just  \ref{item:covariant_representation1}.
For any  $e\in E(S)$  we have $v_e = \pi_e(1_e) = B_*\mhyphen\lim\pi(\mu_i^e)$, where $\{\mu_{i}^e\}_{i}$ is an approximate unit in $I_{e}$.
This $B_*$-convergence written explicitly means that   $f(v_e\xi) = \lim_{i} f(\pi(\mu_i^e)\xi)$  for all $\xi \in E$ and $f\in E'$. This implies that $v_e\xi$, for any $\xi\in E$, sits in the weak closure of $\pi(I_e)E$. As the weak and norm closures of convex sets coincide this means that the range of $v_e$ is  contained in $\overline{\pi(I_{e})E}$.
The reverse inclusion is clear, and so $v_eE= \overline{\pi(I_{e})E}$ for every $e\in E(S)$.
This  is equivalent to \ref{enu:covariant_representations_on_Banach_spaces2}, see Remark \ref{rem:covariant_representations_on_Banach_spaces}.

Conversely, assume that $(\pi,v)$ is a covariant representation on $E$ in the sense of Definition~\ref{defn:covariant_representations_on_Banach_spaces}.
It is clearly a covariant representation in $B(E)$   in the sense of Definition~\ref{defn:covariant_representation_in_Banach_algebras}. Thus we only need to show it is 
normalized.
Let $e\in E(S)$. By Lemma~\ref{lem:covariant_reps_implies_MP^partial_isos}, $\{\pi(\mu_i^e)\}_{i}$ converges strongly to  $v_{e}$, and by construction $\{\pi(\mu_i^e)\}_{i}$ 
$B_*$-converges to $\overline{\pi}_e(1_e)$. Strong convergence is stronger than $B_*$-convergences. Hence $v_{e}=\overline{\pi}_e(1_e)$. 
Using this for any $t\in S$ we get 
$$
v_t= v_{tt^*}v_t= \overline{\pi}_{tt*}(1_{tt*})v_t,
$$
which is equivalent to \eqref{eq:covariant_representation5}.
\end{proof}

\begin{cor}\label{cor:normalization_to_spatial_cov_rep}
If $E$ is reflexive,     for any covariant representation $(\pi,v)$ of $\alpha$ in $B(E)$ there  is a unique covariant representation $(\pi,\widetilde{v})$ on $E$ such that $\pi(a)v_t=\pi(a)\tilde{v}_t$ for all $a\in I_{t}$, $t\in S$.
\end{cor}
\begin{proof}
 By Proposition \ref{prop:reflexive_covariant_representations} and the second part of Corollary \ref{cor:B*_normalized_twisted_rep} the desired pair
 $(\pi,\tilde{v})$ arises as the $E^\prime\widehat{\otimes}E$-normalization of $(\pi,v)$.
\end{proof}
\begin{rem} If the space $E$ has a predual, but is not reflexive, then $B$ has a predual Banach space $B_*$, but is not necessarily a dual Banach algebra.
Then Corollary \ref{cor:predual_Banach_space} applies to normalized covariant representations in $B(E)$, but their relationship  with covariant representations on $E$ is not clear.
\end{rem}
\section{Crossed products for inverse semigroup actions}

Let us fix an inverse semigroup action  $\alpha:S\curvearrowright A$ on a Banach algebra $A$. Following \cite{Sieben97}, we consider the $\ell^1$-direct sum of the ideals $\{I_{t}\}_{t\in S}$ which is the Banach space  
$$
\ell^1(\alpha):=\{f\in\ell^1(S,A):\, f(t)\in I_{t},\,\, t\in S \}
$$
with the norm $\|f\|_1=\sum_{t\in S} \|f(t)\|$. We define multiplication in $\ell^1(\alpha)$ by
$$
(f*g)(r):=\sum_{st=r}\alpha_s(\alpha_{s^*}(f(s))g(t)), \qquad r\in S.
$$ 
It is well defined as for $f,g\in \ell^1(\alpha)$  we have 
\begin{multline*}
\|f*g\|_1=\sum_{r}\|\sum_{st=r}\alpha_s(\alpha_{s^*}(f(s))g(t))\|
\\
\leq \sum_{r}\sum_{st=r}\|f(s)\|\|g(t))\|=\sum_{s,t}\|f(s)\|\|g(t))\|=\|f\|_1\|g\|_1. 
\end{multline*}
The proof of \cite[Proposition 4.1]{Sieben97} shows that this multiplication is associative.
It exploits our assumption that each of the ideals $I_t$, $t\in S$, has an approximate unit. 
In fact without this assumption, associativity of this partial convolution product may fail, see \cite{D-E} for more details on that issue.
This is one of the reasons for including existence of approximate units in the definition of the inverse semigroup action. 

 Hence $\ell^1(\alpha)$ is a Banach algebra. Let $(\pi,v)$ be a covariant representation of $\alpha$ in a Banach algebra $B$. 
For any $f\in\ell^1(\alpha)$ the series 
$$
\pi\times v (f):=\sum_{t\in S}\pi(f(t))v_t
$$
is norm convergent in $B$ and $\|\pi\times v\|\leq \sum_{t\in S}\|\pi(f(t))v_t\|\leq \|f\|_1$. 
Hence $\pi\times v:\ell^1(\alpha)\to B$ is a contractive linear map. 
Using Lemma \ref{lem:range_of_covariant_rep}\ref{enu:range_of_covariant_rep1}
one readily sees that $\pi\times v$ is also multiplicative. Hence $\pi\times v$ 
is a representation. Clearly,  $\overline{\pi\times v (\ell^1(\alpha))} = B(\pi,v)$.
Note that typically $\pi\times v$ will not be injective. 
Indeed, denoting by $a\delta_t$, for any $a\in I_t$, $t\in S$ the element in $\ell^{1}(\alpha)$ such that 
$a\delta_t(s)=0$ for $s\neq t$ and $a\delta_t(t)=a$ we get that $\textup{span} \{a\delta_t: a\in I_{t}, t\in S\}$ is a
dense subalgebra in $\ell^1(\alpha)$. Let us denote by 
$$
\Null:=\clsp\{f(a\delta_s -  a \delta_t)g: f,g\in \ell^{1}(\alpha), a\in I_s, s,t\in S, s\leq t\}
$$ the ideal in $\ell^1(\alpha)$ generated by 
differences $a\delta_s -  a \delta_t$ for $a\in I_s$  and $s,t\in S$ with $s\leq t$.
It follows from Lemma \ref{lem:range_of_covariant_rep}\ref{enu:range_of_covariant_rep2} 
that for any covariant representation $(\pi, v)$ we have 
$$
\Null\subseteq \ker(\pi\times v).
$$
So $\pi\times v$ factors to a representation of the quotient Banach algebra $\ell^{1}(\alpha)/\Null$.
\begin{defn}
    Let $\alpha:S\curvearrowright A$ be an inverse semigroup action. We denote by $A\rtimes_\alpha S$ the Hausdorff completion of $\ell^{1}(\alpha)$ 
			with respect to the submultiplicative seminorm
			$$
			  \|f\|_{\max}:=\sup\{\|\pi\times v(f)\|: (\pi,v) \text{ is a covariant representation of }\alpha \}.
			$$
		We call the arising Banach algebra $A\rtimes_\alpha S$  the \emph{(universal) Banach algebra crossed product} of $\alpha$.
		More generally, if $\mathcal{R}$ is a class of some covariant representations of $\alpha$,  we define the \emph{$\mathcal{R}$-crossed product} of $\alpha$, denoted $A\rtimes_{\alpha,\mathcal{R}} S$, as the Hausdorff completion of $\ell^1(\alpha)$ in the seminorm
    $$
    \|f\|_\mathcal{R}=\sup\{\|\pi\times v(f)\|: (\pi,v)\in\mathcal{R} \}.
    $$
		In particular, any  $(\pi,v)\in\mathcal{R}$  \emph{integrates} to a  representation $\pi\rtimes v$ of $A\rtimes_{\alpha,\mathcal{R}} S$, which sends the image of $a\delta_t$ to  $\pi(a)v_t$, for $a\in I_t$, $t\in S$. 
\end{defn}
\begin{rem}
As we noted above, instead of completions of $\ell^{1}(\alpha)$ one could consider completions of the quotient algebra  $\ell^{1}(\alpha)/\Null$.
In fact, it is often the case (though we do not know in general) that the seminorm $\|\cdot\|_{\max}$ on $\ell^{1}(\alpha)$ becomes a norm on $\ell^{1}(\alpha)/\Null$.
It might even be that $A\rtimes_\alpha S=\ell^{1}(\alpha)/\Null$.
\end{rem}

We now prove the main result of the paper which  gives a universal description of $A\rtimes_\alpha S$ as well `disintegration' of representations of $A\rtimes_\alpha S$
in terms of relevant covariant representations. This result is not present in \cite{Bardadyn-Kwasniewski-McKee1}.
\begin{thm}\label{thm:universal_description_of_crossed_product}
Let $\alpha:S\curvearrowright A$ be an inverse semigroup action on a Banach algebra $A$. The crossed product $A\rtimes_\alpha S$ has the following properties

\begin{enumerate}
\item\label{enu:universal_description_of_crossed_product0} there is a covariant representation 
$(\iota, u)$ of $\alpha$ in $A\rtimes_\alpha S$ (which we may assume to be  normalized) such that $A\rtimes_\alpha S=B(\iota, u)$
\item\label{enu:universal_description_of_crossed_product1} for any covariant representation $(\pi,v)$ in a Banach algebra $B$ there is a unique 
representation $\pi\rtimes v:A\rtimes_\alpha S\to B$ such that 
$$
\pi\rtimes v(\iota(a)u_t)=\pi(a)v_t\quad \text{ for all }a\in I_t,\, t\in S.
$$
\end{enumerate}
Every Banach algebra which satisfying \ref{enu:universal_description_of_crossed_product0} and \ref{enu:universal_description_of_crossed_product1}  above,
is canonically isometrically isomorphic to $A\rtimes_\alpha S$. In addition, 
\begin{enumerate}[label=(\alph*)]
\item\label{enu:universal_description_of_crossed_product2} we have a one-to-one correspondence between representations $\Psi:A\rtimes_\alpha S\to B$ and  normalized covariant representations $(\pi,v)$
of $\alpha$ in $B$  where 
\begin{equation}\label{eq:disintegration_relation}
\Psi(\iota(a)u_t)= \pi(a)v_t  \quad \text{ for all }a\in I_t,\, t\in S.
\end{equation}
If $(B,B_*)$ is a dual Banach algebra, the description of normalized covariant representations of $\alpha$ in $B$ simplifies to the one in Corollary \ref{cor:B*_normalized_twisted_rep} 
\item\label{enu:universal_description_of_crossed_product4} If $E$ is a reflexive Banach space, then \eqref{eq:disintegration_relation} establishes a one-to-one correspondence between representations $\Psi:A\rtimes_\alpha S\to B(E)$ and  covariant representation $(\pi,v)$
of $\alpha$ on $E$. 
\end{enumerate}
\end{thm}
\begin{proof}

By construction, for each $a\in \ell^{1}(\alpha)$ and $\varepsilon >0$ there is a covariant representation $(\pi,v)$ such that $\|a\|_{\max}<\|\pi\times v(a)\|+\varepsilon$.
Thus there exist a family $\{(\pi^{i},v^{i})\}_{i\in I}$ where $(\pi^{i},v^{i})$ is a covariant representation of $\alpha$ in a Banach algebra  $B_i$ 
and such that $\|a\|_{\max}=\sup_{i\in I}\|\pi\times v(a)\|$ for every $a\in \ell^{1}(\alpha)$. 
Let us now express it using  the $\ell^{\infty}$-sum of these representations. Namely, let 
$$
B:=\oplus^{\ell^{\infty}}_{i\in I} B_i''=\{f:I\to \bigsqcup_{i} B_i'': f(i)\in B_{i}'', i\in I, \|f\|_{\infty}=\sup\|f(i)\|<\infty\}.
$$
We define the algebra homomorphism $\iota:A\to B$ by $\iota(a)=\oplus_{i\in I}^{\ell^{\infty}}\pi_i(a)$, $a\in A$, and the 
semigroup homomorphism $u:S\to B$ by $u_t:=\oplus_{i\in I}^{\ell^{\infty}}\ v^{i}_t$, $t\in S$. 
Clearly, $(\iota, u)$ is a covariant representation of $\alpha$ in $B$ such that 
$\iota\times  u=\oplus_{i\in I}^{\ell^{\infty}}\pi^{i}\times v^{i}_t$ and so $\|\iota\times  u(a)\|=\|a\|_{\max}$ for $a\in\ell^{1}(\alpha)$. 
Thus it induces an isometric representation $\iota\rtimes  u$ of $A\rtimes_\alpha S$ and so we may identify its range with $A\rtimes_\alpha S$.
Then $(\iota, u)$ is a covariant representation  in $A\rtimes_\alpha S$ such that $A\rtimes_\alpha S=B(\iota, u)$. 
Property \ref{enu:universal_description_of_crossed_product1} is clear by construction. 
If $C=B(\tilde{\iota}, \tilde{u})$ is the range of a  covariant representation $(\tilde{\iota}, \tilde{u})$ satisfying an analogue of \ref{enu:universal_description_of_crossed_product1}
then $\iota\rtimes u$ and  $\tilde{\iota}\rtimes  \tilde{u}$ are mutually inverse representations. Hence they are isometric and $A\rtimes_\alpha S\cong C$.

For \ref{enu:universal_description_of_crossed_product2}, let
$\Psi:A\rtimes_\alpha S\to B$ be a representation, and consider the representation $\Psi'':(A\rtimes_\alpha S)''\to B''$.
Then $(\Psi''\circ \iota, \Psi''\circ u)$ is a well defined covariant representation of $\alpha$ in $B$
such that $\Psi= (\Psi''\circ \iota)\rtimes (\Psi''\circ u)$.
By passing to  $B'$-normalization, see Proposition 
\ref{prop:normalized_twisted_rep},  we get \ref{enu:universal_description_of_crossed_product2}.
By Proposition \ref{prop:reflexive_covariant_representations}, \ref{enu:universal_description_of_crossed_product4} is a special case of \ref{enu:universal_description_of_crossed_product2}.
\end{proof}

\begin{cor}\label{cor:group_actions}
Let $\alpha: G\curvearrowright A$ be an action of a group $G$  (so it is also an inverse semigroup action).
Then  the group crossed product
$\ell^1(G,A)$, the algebra $\ell^{1}(\alpha)$ and the inverse semigroup crossed product $A\rtimes_\alpha G$ coincide. 

Moreover, covariant representations of $\alpha$ as defined on page  \pageref{page:covariant_rep} coincide with non-degenerate normalized covariant representations of  $\alpha$ treated as the inverse semigroup action.
\end{cor}
\begin{proof}
The first part is straightforward. 
For the second part it suffices to show that if $(\pi,v)$ is  a  non-degenerate normalized covariant representation $(\pi,v)$ of $\alpha$ in $B$, then we necessarily have that the semigroup homomorphism $v:G\to  (B'')_1$ takes values in the group $\UMult(B)$.
But by  definition of the covariant representation we have that $v_{g}\pi(a)=\pi(\alpha_g(a))v_g\in B$ and $\pi(a)v_{g}\in B$  for $g\in G$ and $a\in A$.
Since $\pi(A)B=B=B\pi(A)$ by non-degeneracy of $\pi$, this implies that $v_gB\subseteq B$ and $Bv_g\subseteq B$. Hence  $\{v_g\}_{g\in G}\subseteq  \Mult(B)$.
By $B'$-normalization $v_1$ is the unit in by $\Mult(B)$. Thus $\{v_g\}_{g\in G}\subseteq  \UMult(B)$ as they are contractive and invertible.
\end{proof}

Finally let us comment on $C^*$-algebraic crossed products. Let $A$ be a $C^*$-algebra.  
Then the Banach algebra $\ell^1(\alpha)$ is naturally a Banach $*$-algebra with involution
given by 
$$
f^*(t) := \alpha_{t^*}(f(t^*)^*),\qquad  t\in S.
$$
If $(\pi,v)$ is a covariant representation of $\alpha$ in a $C^*$-algebra $B$, then Corollary \ref{cor:range_of_covariant_rep} implies that the representation $\pi\times v:\ell^{1}(\alpha)\to B$ is 
automatically a $*$-homomorphism.
Sieben \cite[Definition 4.4]{Sieben97} defined the $C^*$-algebraic crossed product   by 
completing $\ell^{1}(\alpha)$ using $*$-homomorphism $\pi\times v$ coming from  covariant representations $(\pi,v)$ on Hilbert spaces, as we described in Definition \ref{defn:covariant_representations_on_Banach_spaces}.
 Buss and Exel \cite[Page 257]{Buss-Exel} defined the $C^*$-algebraic crossed product  using covariant representations $(\pi,v)$ in $C^*$-algebras satisfying properties in Proposition \ref{prop:normalized_twisted_rep}, which characterize what we call normalized covariant representations. Using our normalization procedure we can now explain why these two  different constructions coincide. 
\begin{cor}\label{cor:Sieben_Buss}
    Let $\alpha: S\curvearrowright A$ be an action of an inverse semigroup $S$ where $A$ is a $C^*$-algebra. Let $\mathcal{R}_{C^*}$ be the class of all covariant representations $\alpha$ in some $C^*$-algebra and $\mathcal{R}_{\text{H}}$ be the class of all covariant representations of $\alpha$ on some Hilbert space. Then
    $$
    A\rtimes_{\alpha,\mathcal{R}_{C^*}}S=A\rtimes_{\alpha,\mathcal{R}_{\text{H}}}S.
    $$
This algebra is a  $C^*$-algebra which coincides with the crossed products introduced in \cite{Sieben97} and \cite{Buss-Exel}.
\end{cor}
\begin{proof}
Since $\RR_{\textup{H}} \subseteq \RR_{C^*}$ we get $\|\cdot\|_{\RR_{\textup{H}}} \leq \|\cdot\|_{\RR_{C^*}}$. 
This is equality, because for any covariant representation $(\pi,v)\in \RR_{C^*}$ in a $C^*$-algebra $B$, we may assume that the enveloping $W^*$-algebra $B''$ is faithfully represented on a Hilbert space $H$. 
Then by Corollary~\ref{cor:normalization_to_spatial_cov_rep} the normalized covariant representation  $(\pi,\tilde{v})$ is a  covariant representation of $\alpha$ on $H$, so $(\pi,\tilde{v})\in \RR_{\textup{H}}$, and   $\pi\times v(f) = \pi\times \tilde{v}(f)$ for all $f \in \ell^1(\alpha)$.
Hence $\|\cdot\|_{\RR_{\textup{H}}}=\|\cdot\|_{\RR_{C^*}}$. 
 By Corollary \ref{cor:range_of_covariant_rep}, this is a $C^*$-seminorm on  $\ell^{1}(\alpha)$, which coincides with those used in \cite{Sieben97} and \cite{Buss-Exel}.
\end{proof}

\subsection*{Acknowledgment}
We thank the anonymous referee for their comments. 
The second named author  was  supported by the National Science Centre, Poland, through the WEAVE-UNISONO grant no. 2023/05/Y/ST1/00046. 


\bibliographystyle{spmpsci}
\bibliography{myBibLib} 

\end{document}